\documentclass[a4paper,11pt]{article}
\usepackage{amsmath,amscd,amsfonts,amssymb,epsf,latexsym}

\makeatletter

\long\def\@makefnt#1{\parindent 1em\noindent
             \hb@xt@1.8em{\hss\@textsuperscript{}}#1}
\long\def\@ftntext#1{\insert\footins{%
     \reset@font\footnotesize
     \interlinepenalty\interfootnotelinepenalty
     \splittopskip\footnotesep
     \splitmaxdepth \dp\strutbox \floatingpenalty \@MM
     \hsize\columnwidth \@parboxrestore
     \color@begingroup
       \@makefnt{%
         \rule\z@\footnotesep\ignorespaces#1\@finalstrut\strutbox}%
     \color@endgroup}}%
\def\subjclass#1{%
   \@ftntext{2000 {\itshape Mathematics Subject Classification.}\enspace
#1.}}
\def\keywords#1{%
   \@ftntext{{\itshape Key words and phrases.}\enspace #1.}}
\makeatother

\def\moins{\raise 1pt\hbox{{$\scriptstyle -$}}}
\def\plus{\raise 1pt\hbox{{$\scriptstyle +$}} }

\newtheorem{theorem}{Theorem}
\newtheorem{proposition}[theorem]{Proposition}
\newtheorem{lemma}[theorem]{Lemma}
\newtheorem{corollary}[theorem]{Corollary}
\newtheorem{remark}[theorem]{Remark}

\newtheorem{procedure}[theorem]{Procedure}
\newtheorem{example}[theorem]{Example}

\def\proof{\noindent{\bf Proof.\ }}

\def\qed{~\hbox{$\Box$}}

\def\Sym{\mathop{\rm Sym}}

\def\cT{{\mathcal T}}

\def\cO{{\mathcal O}}
\def\C{{\mathbb C}}
\def\Z{{\mathbb Z}}
\def\N{{\mathbb N}}
\def\Qm{{\mathbb Q}}
\def\Q{\widetilde{Q}}
\def\gJ{{\frak J}}

\def\s{\vskip6pt}

\def\P{{\bf P}}
\def\qed{~\hbox{$\Box$}}

\begin{document}

\title{\bf Positivity of Legendrian Thom polynomials}

\author{Ma\l gorzata Mikosz\thanks{Research supported by the KBN grant NN201 387034.}\\
\small Warsaw University of Technology,\\
\small Pl.~Politechniki 1, 00-661, Warszawa, Poland\\
\small emmikosz@mech.pw.edu.pl \and
Piotr Pragacz\thanks{Research supported by the MNiSzW grant N N201 608040.}\\
\small Institute of Mathematics of Polish Academy of Sciences\\
\small \'Sniadeckich 8, 00-956 Warszawa, Poland\\
\small P.Pragacz@impan.pl \and Andrzej Weber\thanks{Research
supported by the KBN grant NN201 387034.}\\
\small Institute of Mathematics, University of Warsaw\\
\small Banacha 2, 02-097 Warszawa, Poland\\
\small aweber@mimuw.edu.pl \and
\footnotesize and: IMPAN, \'Sniadeckich 8, 00-956 Warszawa, Poland}

\subjclass{05E05, 14C17, 14N15, 55R40, 57R45}

\keywords{Legendre and Lagrange singularities, Thom polynomials, $\Q$-functions,
Schur $Q$-functions, jet bundles, numerical positivity}


\maketitle


\begin{abstract} We study Legendrian singularities arising in
complex contact geometry. We define a one-parameter family of
bases in the ring of Legendrian characteristic classes such that
any Legendrian Thom polynomial has nonnegative coefficients when
expanded in these bases. The method uses a suitable Lagrange
Grassmann bundle on the product of projective spaces. This is an
extension of a nonnegativity result for Lagrangian Thom
polynomials obtained by the authors previously. For a fixed
specialization, other specializations of the parameter lead to
upper bounds for the coefficients of the given basis. One gets
also upper bounds of the coefficients from the positivity of
classical Thom polynomials (for mappings), obtained previously by
the last two authors.

\end{abstract}

\section{Introduction}

The aim of the present paper (which is a continuation of
\cite{Ka2}, \cite{PW} and \cite{MPW}) is to study the {\it positivity}
of Legendrian Thom polynomials. The pioneering papers \cite{G}
of Griffiths and \cite{FL} of Fulton and Lazarsfeld investigated
numerical positivity related to ample vector bundles in differential
and algebraic geometry, respectively. Their various
variants are nowadays widely investigated in algebraic geometry.
We refer to the monograph \cite{Lz} for more detailed account.
Also, the recent paper \cite{FNR} gives a combinatorial interpretation
for the coefficients of certain Thom-like polynomials, providing
further motivation for studying general positivity phenomena.

Our main goal is to define a certain one-parameter family of bases in the ring of
Legendrian characteristic classes. Any Legendrian Thom polynomial
has nonnegative  coefficients when expanded in any member of this
family (Theorems \ref{main2} and \ref{mainfamily}).

The main difference, comparing with the previous papers (in particular,
with \cite{Ka2} and \cite[Remark 14]{MPW}) is the definition of
{\it Legendre singularity classes} to which these Thom polynomials are attached.
They are introduced as closed algebraic subsets in the space parametrizing
{\it pairs} of Legendrian submanifolds in the contact space (see Lemma \ref{LL}).
Regardless of the technical differences the main theorem of \cite{MPW} follows from
Theorem \ref{mainfamily} by specialization of the parameter $t$ to $0$.
As in \cite{Ka2}, \cite{PW} and \cite{MPW}, when we study Thom polynomials of commonly
bounded degree, it is enough to deal with finite jets of germs of submanifolds through the origin.

The principal technique involved in the proof of these theorems is {\it trans\-versality}
with respect to some stratification of a Lagrange Grassmann bundle.
This is a subject of Section \ref{pos}.
The key technical result is Theorem \ref{trans1}. We show that the intersection
(in the jet bundle) of a Legendre singularity class with the preimage of the closue
of a stratum of the stratification, is represented by a nonnegative cycle.
In fact, for our purposes we need a certain Lagrange Grassmann bundle over
the product of two projective $n$-spaces, see Theorem \ref{main2}.
A concrete form of the resulting one-parameter family of bases follows from degeneracy
locus formulas from \cite{KT}, \cite{LP}, \cite{PR} and a formula of Kazarian.

In Section \ref{parameter}, we examine the parameter of the constructed
family of bases, and give a precise proof of a result announced in \cite{MPW},
describing algebraically the basis corresponding to the value of the parameter
equal to 1 (see Theorem \ref{Qt}).

Let us fix the value of the parameter. It turns out that the nonnegativity
of coefficients of the bases for some other values of the parameter
can imply upper bounds of the coefficients in the basis for the given parameter.
This is a subject of Section \ref{bounds}.

In Section \ref{vs}, we show (Proposition \ref{nonzero}) that for
nonempty stable singularity classes the corresponding Legendrian
(and Lagrangian) Thom polynomials are nonzero. This is an
amelioration of the main result of \cite{MPW}. The proof uses the
fact \cite{PW} that the Thom polynomials for functions ${\Bbb C}^n
\to {\Bbb C}$ are nonzero for nonnempty singularity classes. We
also show how this last result gives some upper bounds on the
coefficients of Legendrian Thom polynomials in the basis from
Theorem \ref{Qt}.

In Section \ref{examples}, we list some examples of Legendrian Thom
polynomials expanded in different bases from the family.

In the Appendix, we prove some new positivity result for
the intersection coefficients in the Lagrange Grassmann
bundle. This result concerns the stratification transverse
to the one studied in Section \ref{pos}. It relies
on a``large group'' action from \cite{An}.

\section{Some Legendrian geometry}

Fix $n \in \Bbb N$.
Suppose for the moment that  $W$ is a vector space of dimension $n$, and $\xi$ is a
vector space of dimension one. Let
\begin{equation}
V:=W\oplus (W^*\otimes\xi)
\end{equation}
be the standard symplectic space equipped with the twisted symplectic form
$\omega\in\Lambda^2V^*\otimes\xi$.
We study the germs at the origin of the Legendrian submanifolds in the standard {\it contact space}
$$
V\oplus\xi=W\oplus (W^*\otimes\xi)\oplus \xi\,,
$$
or equivalently the germs of the Lagrangian submanifolds in the symplectic space $V$.
Any Legendrian submanifold in $V\oplus\xi$ is determined by its projection to $V$ and any Lagrangian submanifold in $V$ lifts to $V\oplus\xi$, see \cite[Proposition p.313]{AGV}. Therefore we will perform all the constructions in the realm of symplectic geometry.

Two Lagrangian submanifolds, if they are in generic position, intersect transversally.  The singular  relative positions can be divided into
Legendrian singularity classes. To classify all the possible relative positions it suffices consider only two types of submanifolds:
\begin{description}
\item {(i)} {\it linear subspaces}; they are parametrized by Lagrangian Grassmannian denoted by $LG(V,\omega)$;
\item {(ii)} {\it the submanifolds which have the tangent space at the origin equal to $W$}; they are the graphs of the differentials of the functions $f:W\to\xi$ satisfying $df(0)=0$ and $d^2f(0)=0$, see \cite[Lemma 2]{MPW}.
\end{description}

The group of symplectomorphisms of $V$ acts on the pairs of Lagrangian subspaces.
By the following lemma we can restrict our attention only to the submanifolds of the type (i) and (ii).

\begin{lemma}\label{LL}
Any pair of Lagrangian submanifolds is symplectic equivalent to a pair $(L_1,L_2)$ such
that $L_1$ is a linear Lagrangian subspace and the tangent space $T_0L_2$ is equal to $W$.
\end{lemma}

\medskip

\proof The lemma follows essentially from the Darboux Theorem (see \cite{AGV}, Theorem, p.287).
Indeed, it follows from this theorem that any Lagrangian submanifold is symplectomorphic
to a linear one (given by vanishing of the $p$-coordinates in the notation of the theorem).
We apply this to the first Lagrangian submanifold, getting a linear $L_1'$.
Applying then an appropriate rotation, we get the tangency condition for $L_2$, equal to
the rotated second Lagrangian submanifold. The image $L_1$ of $L_1'$ under this rotation is linear.
The pair $(L_1,L_2)$ satisfies the assertion of the lemma.
\qed

\medskip

Let us fix a suitable large $k$.
 We identify two Lagrangian submanifolds if the degree of their
tangency at $0$ is greater than  $k$. The equivalence class will be called
``a $k$-jet of a submanifold''. The $k$-jets of submanifolds satisfying the condition (ii) are in bijection with elements of the vector space
 $$\bigoplus_{i=3}^{k+1}{\Sym}^i(W^*)\otimes\xi\,.$$

\medskip

We want to describe a space parametrizing all possible relative positions of Lagrangian submanifolds.
A suggestion of Kazarian and Lemma \ref{LL} justifies the following definition.
We denote by $\gJ^k(W,\xi)$ the set of pairs $(L_1,L_2)$ of $k$-jets of Lagrangian submanifolds of $V$
such that  $L_1$ is a linear space and $T_0L_2=W$. Let
\begin{equation}\label{pi}
\pi: {\gJ}^k(W,\xi) \to LG(V,\omega)
\end{equation}
be the projection to the first factor. Clearly, $\pi$ is a trivial vector bundle with the fiber equal to:
$$
\bigoplus_{i=3}^{k+1}{\Sym}^i(W^*)\otimes\xi \,.
$$

\begin{remark}\rm In \cite{MPW} we considered the relative position of the Lagrangian submanifolds with respect to the fixed linear space $W^*$. We obtained a jet bundle over $LG(V)$ which was not a trivial bundle, and we had to deform that bundle to its linear part.
Fixing the tangent space of the nonlinear submanifold and moving the second linear space, is an important simplification comparing with \cite{MPW}.
\end{remark}

\medskip
We are interested in a larger group than just the group of symplectomorphisms, namely
the group of {\it (complex) contact automorphisms} of $V\oplus \xi$. It acts on the pairs of
Legendrian submanifolds in $V\oplus \xi$. Again by \cite[Proposition p.313]{AGV}, we obtain an action on the pairs of Lagrangian submanifolds. In particular, we take into  account the automorphisms of $V$ which transport the symplectic form $\omega$ to a proportional one. For example,
$$
\phi_t(q,p)=(q,tp), \ \ \ \hbox{where} \ \ \ q\in W, \ p\in W^*\otimes \xi \ \ \hbox{and} \ \ t\in \C^*.
$$
By a {\it Legendre singularity class} we mean a closed algebraic subset
$$
\Sigma \subset {\gJ}^k(\C^n,\C)\,,
$$
invariant with respect to holomorphic contactomorphisms of $\C^{2n+1}$.
 Additionally, we assume that
the singularity class $\Sigma$ is stable with respect to enlarging the
dimension of $W$, as  in \cite[Sect.2]{PW}.
Unfortunately, we do not know any place in literature where the relation between cohomological stability and the infinitesimal stability in the sense of see \cite[Sect.6]{AGV} is discussed. This problem is treated in \cite[Sect. 7.2]{FR} in a different context.
It seems to be a common knowledge that infinitesimal stability implies cohomological stability.
On the other hand, we would like to mention that our main results: Theorems \ref{main2} and \ref{mainfamily} hold without the stability assumption.

\section{The jet bundle $\gJ^k(W,\xi)$}

The vector space $\xi$ may have no distinguished coordinate. It
happens so for example  when we deal with a fiber of a vector
bundle. In other words we have a nontrivial action of $\C^*$ on
${\gJ}^k(\C^n,\C)$. Now we repeat the construction of 
the space
$
{\gJ}^k(W,\xi)
$
parametrizing the relative positions of two Lagrangian submanifolds,
assuming that $\xi$ is a line bundle over some
base space. We could have used the universal base $BU(1)$  but it
is more convenient to work with bundles defined over various base
spaces. Also, it will be useful to assume that $W$ is a (possibly
nontrivial) vector bundle.

Let $X$ be a topological space, $W$ a complex rank $n$ vector bundle over $X$, and $\xi$ a complex
line bundle over $X$. The fibers of $W$, $\xi$, $V$ ~over a point $x\in X$ are denoted
by $W_x$, $\xi_x$, $V_x$. Let
\begin{equation}
\tau: LG(V,\omega)\to X
\end{equation}
denote the Lagrange Grassmann bundle parametrizing Lagrangian
linear submanifolds in $V_x$, $x\in X$. We have a relative version of the map (\ref{pi}):
\begin{equation}
\pi: {\gJ}^k(W,\xi) \to LG(V,\omega)
\end{equation} 
The space ${\gJ}^k(W,\xi)$ fibers over $X$. It is equal to the pull-back:
\begin{equation}
{\gJ}^k(W,\xi) = \tau^*\left(\bigoplus_{i=3}^{k+1}{\Sym}^i(W^*)\otimes\xi\right)\,.
\end{equation}
In the following,
the pull-backs to $LG(V,\omega)$ of the bundles $W$, $V$ and $\xi$ will be denoted by the same
letters, if no confusion occurs.
Since any changes of coordinates of $W$ and $\xi$ induce
holomorphic contactomorphisms of $V\oplus \xi$, any Legendre
singularity class $\Sigma$ defines a cycle
\begin{equation}
\Sigma(W,\xi)\subset{\gJ}^k(W,\xi).
\end{equation}
We will study the classes\footnote{In this paper, whenever we speak about the classes of algebraic cycles, we always
mean their {\it Poincar\'e dual classes} in cohomology.} defined by the cycles $\Sigma(W,\xi)$.

\section{Legendrian characteristic classes}

The tautological bundle over $LG(V,\omega)$ is denoted by $R_{W,\xi}$, or by $R$ for short. The symplectic form $\omega$ gives an isomorphism
\begin{equation}
V\cong V^*\otimes\xi\,.
\end{equation}
There is a tautological sequence of vector bundles on $LG(V,\omega)$:
\begin{equation}\label{taut} 0\to R \to V\to R^*\otimes\xi\to 0\,.
\end{equation}
 Consider the virtual bundle
\begin{equation}\label{A}
A:=W^*\otimes \xi-R_{W,\xi}\,.
\end{equation}
Using the sequence (\ref{taut}), we get the relation
\begin{equation}\label{rel}
A+ A^*\otimes \xi=0\,.
\end{equation}
The Chern classes  $a_i=c_i(A)$  generate the cohomology
$H^*(LG(V,\omega),\Z)\cong H^*({\gJ}^k(W,\xi), \Z)$ as an algebra
over $H^*(X,\Z)$.

\medskip

Let us fix an approximation of $BU(1)=\bigcup_{n\in\N}\P^n$, that is we set $X=\P^n$, $\xi=\cO(1)$.
Let $W={\bf 1}^n$ be the trivial bundle of rank $n$.
Then the cohomology
$$
H^*(LG(V,\omega), \Z)\cong H^*({\gJ}^k(W,\xi),\Z)
$$
is isomorphic  to the ring of  Legendrian characteristic classes
for degrees smaller than or equal to $n$. The element
$[\Sigma(W,\xi)]$ of $H^*({\gJ}^k(W,\xi), \Z)$,
is called the {\it Legendrian Thom
polynomial} of $\Sigma$, and is often denoted by ${\cal T}^{\Sigma}$. It is written in terms of the generators
$a_i$ and $s=c_1(\xi)$,
 \ (cf.
\cite[Sect. 3.4]{Ka2}, \cite[Sect.4]{Kahab}).

\begin{remark}\label{Legring} \rm The ring of Legendrian characteristic
classes is the quotient of the polynomial ring
$$
\Z[a_1,a_2,a_3,\dots;u]
$$
by the relations coming from the identity (\ref{rel}). After inverting
2 and applying the twist by $\xi^{-\frac12}$ we obtain the ring
of Lagrangian classes extended by one additional free variable $t$,
that is
$$
\Z[\frac12][a'_1,a'_2,a'_3,\dots;t]/{\cal I}\,,
$$
where ${\cal I}$ is generated by the polynomials
$$
(a'_i)^2+2\sum_{k=1}^i(-1)^ka'_{i+k}a'_{i-k},\qquad i>0\,.
$$
The even Chern classes $a'_{2i}$ are expressed by odd ones and
this ring is just the polynomial ring
$$
\Z[\frac12][a'_1,a'_3,a'_5,\dots;t]\,.
$$
A similar procedure can be applied to the untwisted variables
$a_i$.
\end{remark}

\section{Cell decompositions of the  Grassmann bundle}\label{LegGrass}

We describe two ``transverse'' cell decompositions of the Lagrange
Grassmannians.

To begin with, let $\xi,\alpha_1,\alpha_2,\dots,\alpha_n$ be
vector spaces of dimensions equal to one and let
\begin{equation}
W:=\bigoplus_{i=1}^n\alpha_i\,,\qquad V:=W\oplus (W^*\otimes \xi)\,.
\end{equation}
 We have a twisted symplectic form $\omega$ defined on $V$ with
values in $\xi$. The Lagrangian Grassmannian $LG(V,\omega)$
is a homogeneous space with respect to the group action of
the symplectic group $Sp(V,\omega)\subset {\rm End}(V)$.
Fix two ``opposite'' standard isotropic flags in $V$:
\begin{equation}
F_h^+:=\bigoplus_{i=1}^h\alpha_i\,,\qquad
F_h^-:=\bigoplus_{i=1}^h\alpha_{n-i+1}^*\otimes\xi
\end{equation}
for $h=1,2,\dots ,n$ and consider two subgroups $B^\pm\subset
Sp(V,\omega)$ which are the Borel groups preserving the flags
$F^\pm_{\bullet}$. The orbits of $B^\pm$ in $ LG(V,\omega)$ form
two ``opposite'' cell $\Omega^{\pm}$-decompositions $\{\Omega_I(F^\pm_\bullet,\xi)\}$
of $LG(V,\omega)$ (Bruhat decompositions), indexed by strict partitions
$$
I\subset \rho:=(n,n-1,\ldots,1)
$$
(see, e.g., \cite{P1}). The cells of the $\Omega^+$-decomposition are
transverse to the cells of the $\Omega^-$-decomposition.

We pass now to the relative version of the above decompositions.

The description just presented is functorial with respect to the
automorphisms of the lines $\xi$ and $\alpha_i$'s,
(they form a torus $(\C^*)^{n+1}$). Thus the construction of the cell
decompositions can be repeated for bundles $\xi$ and
$\{\alpha_i\}_{i=1}^n$ over any base $X$. We obtain a Lagrange Grassmann
bundle
$$
\tau: LG(V,\omega)\to X
$$
and a group bundle (group scheme over $X$)
$$
Sp(V,\omega)\to X
$$
together with two subgroup bundles $B^\pm \to X$. Moreover, $LG(V,\omega)$
admits two (relative) stratifications
$$
\{\Omega_I(F^\pm_\bullet,\xi) \to X\}_{{\rm strict} \ I\subset \rho}\,.
$$

Assume that $X=G/P$ is a compact manifold, homogeneous with
respect to an action of a linear group $G$. Then $X$
admits an algebraic cell decomposition $\{\sigma_\lambda\}$ (it is
again a Bruhat decomposition). The subsets
\begin{equation}
Z^-_{I\lambda}:=\tau^{-1}(\sigma_\lambda)\cap\Omega_I(F^-_\bullet,\xi)
\end{equation}
form an algebraic cell decomposition of $LG(V,\omega)$, called
$Z^-$-{\it decomposition} or
{\it distinguished decomposition} in the following. The classes of their
closures give a basis of homology, called $Z^-$-basis. Note that each $Z^-_{I\lambda}$
is transverse to each stratum $\Omega_J(F^+_\bullet,\xi)$,
where $J\subset \rho$ is a strict partition.

Similarly, we define the subsets
\begin{equation}
Z^+_{I\lambda}:=\tau^{-1}(\sigma_\lambda)\cap \Omega_I(F^+_\bullet,\xi)\,.
\end{equation}
which form a $Z^-$-{\it decomposition} and give rise to the corresponding
$Z^-$-basis of the cohomology.

\begin{example}\label{Hirz} \rm
If $X=\P^1$, $W={\bf 1}$, $\xi=\cO(d)$ (for $d>0$),
then $LG(V,\omega)$ is the Hirzebruch surface $\Sigma_d$ which can be
presented as the sum of the space of the line bundle $\xi$ and
the section at the infinity,
$$
\Sigma_d=\xi \cup \P^1_\infty\,.
$$
Then $\P^1_0$, the zero section of the bundle $\xi$, is a stratum
of the $\Omega^+$-decomposition and the section at infinity
$\P^1_\infty$ is a stratum of the $\Omega^-$-decomposition. For
the cell decomposition of $X=\P^1=\C\cup \{\infty\}$, we obtain
two cell decompositions of $\Sigma_d$. The closures of the
one-dimensional cells are the following:
$$
  \Big\{\overline {Z^{\,+}_{I\lambda}}\Big\}=
  \Big\{\tau^{-1}(\infty),\;\P^1_0\Big\}\,,
  \qquad
  \Big\{\overline {Z^{\,-}_{I\lambda}}\Big\}=
  \Big\{ \tau^{-1}(\infty),\;\P^1_\infty\Big\}\,.
$$
Two resulting bases of {\it cohomology} are mutually dual with
respect to the intersection product.

The cycles of $Z^+$-decomposition have the following property:
any effective cycle has a nonnegative intersection number with them.
This is not true for the elements of $Z^-$-decomposition: for example, the
self-intersection of $\P^1_\infty$ is equal to $-d$.
\end{example}

\section{Positivity in a jet bundle on a Legendrian Grassmannian}\label{pos}

We pass now to a nonnegativity result on the Legendrian Thom
polynomials and the $Z^{-}$-decomposition. We recall that our
goal is to study cycles $\Sigma(W,\xi)$ in:
$$
{\gJ}^k(W,\xi)=\tau^*\left(\bigoplus_{i=3}^{k+1}{\Sym}^i(W^*)\otimes\xi\right).
$$
To abbreviate the notation, we set ${\gJ}^k(W,\xi)=\gJ$.

This vector bundle is equipped with a
$B^+$-action, i.e. an appropriate map
$$
B^+\times_X \gJ\to \gJ
$$
over $X$.

A Legendre singularity class $\Sigma$ defines the subset
$\Sigma(W,\xi)\subset \gJ$, which is a Zariski-locally trivial
fibration over $X$. Moreover, $\Sigma(W,\xi)$ is preserved by the
action of $B^+$ since this group consists of holomorphic
symplectomorphisms preserving $W$:

$$\begin{matrix}&&&\pi&&\tau\cr
 \Sigma(W,\xi)&\subset&\gJ&\twoheadrightarrow&LG(V,\omega)&
\twoheadrightarrow&X\cr
 \circlearrowleft&&\circlearrowleft&&\circlearrowleft\cr
 B^+&&B^+&&Sp(V,\omega)\cr
 \end{matrix}$$

We now state the main technical result of the present paper.

\begin{theorem}\label{trans1} Fix $I\subset \rho$ and $\lambda$.
Suppose that the vector bundle $\gJ$ is globally generated. Then,
in $\gJ$, the intersection of $\Sigma(W,\xi)$ with the closure of
any $\pi^{-1}({Z^-_{I\lambda}})$ is represented by a nonnegative
cycle.
\end{theorem}
This result is, in fact, true for any effective $B^+$-invariant
cycle on $\gJ$, Zariski-locally trivial fibered over $X$, taken instead of
$\Sigma(W,\xi)$. The proof is modelled on the techniques of
\cite{Kl}.

Before proving the theorem, we shall establish some preliminary
result. For a subset $Y \subset \gJ$ and a global section $s\in
H^0(\gJ)$, we denote by $s+Y$ the fiberwise translate of $Y$ by
$s$. We will deform the cycle $\Sigma(W,\xi)$ using a fiberwise
translate. The construction is done for each stratum
$\Omega_J(F^+_{\bullet},\xi)$ separately. Fix such a stratum
$\Omega=\Omega_J(F^+_{\bullet},\xi)$. Denote by
$\gJ_{|\Omega}=\pi^{-1}\Omega$ the restriction of the bundle to
the stratum and set
\begin{equation}
\Xi:=\Sigma(W,\xi)\cap\pi^{-1}\Omega\,.
\end{equation}

We need the following lemma which is a variant of Bertini-Kleiman tranversality theorem \cite{Kl}.

\begin{lemma}\label{lem} Suppose that the vector bundle $\gJ$ is globally generated. Let $Y\subset \gJ_{|\Omega}$ be a
subvariety.  Then there exists an open, dense subset $U\subset
H^0(\gJ)$ such that for any section $s\in U$, the translate
$s+\Xi$ has proper intersection with $Y$.
\end{lemma}
\proof
Let
$$
q:H^0(\gJ)\times \Xi \to \gJ_{|\Omega}
$$
be the fiberwise translate. We claim that this map is flat (in
fact, it is a fibration). The question is local. We find an open
set $X'\subset X$ such that the bundles ${\alpha_i}_{|X'}$ and
$\xi_{|X'}$ are trivial. Over $X'$ the set $\Xi$ is the product of
$X'$ and some variety. Therefore, it is enough to assume that $X$
is a point.  Further, note that $\Xi\to\Omega$ is a fibration
since $\Omega$ is homogeneous with respect to $B^+$ and $\Xi$ is a
$B^+$-invariant subset of $\gJ$. Thus the question reduces to a
single fiber of $\gJ$. Since $\gJ$ is globally generated, for
any $y\in \Omega$ the fiber $\gJ_y$ is homogeneous for the action
of $H^0(\gJ)$. The action map
$$
q_{|H^0(\gJ)\times\Xi_y}:H^0(\gJ)\times\Xi_y \rightarrow \gJ_y
$$
is  a trivial fibration with the fiber isomorphic to $\Xi_y\times
{\rm ker}(H^0(\gJ)\to\gJ_y)$. It follows that $q$ is a fibration.

Applying \cite[Lemma 1]{Kl}, the assertion of the lemma follows.
\qed

\bigskip

\noindent {\bf Proof of Theorem \ref{trans1}}

\smallskip

For a strict partition $J\subset \rho$, we set
\begin{equation}
\gJ_    J:=\pi^{-1}(\Omega_J(F^+_\bullet,\xi))\,.
\end{equation}
Applying Lemma \ref{lem}, we get an open dense subset $U_J\subset
H^0(\gJ)$ such that for $s\in U_J$ the intersection
$$
\big(s+(\Sigma(W,\xi)\cap
\gJ_J)\big)\cap\big(\pi^{-1}(\overline{Z^-_{I\lambda}})\cap
\gJ_J\big)
$$
is proper inside $\gJ_J$. We now pick
$$
s\in\bigcap_{{\rm strict} \ J \subset \rho} U_J\,,
$$
and set
\begin{equation}
\Sigma(W,\xi)':=s+\Sigma(W,\xi)\,.
\end{equation}
Since the cycle $\pi^{-1}(\overline{Z^-_{I\lambda}})$ is
transverse to the stratification $\{\gJ_J\}_{{\rm strict} \ J
\subset \rho}$ of $\gJ$, an easy dimension count shows that
$\pi^{-1}(\overline{Z^-_{I\lambda}})$ intersects properly
$\Sigma(W,\xi)'$ in $\gJ$.

Theorem \ref{trans1} now follows since by \cite[Sect. 8.2]{Fu} the intersection
$$
[\Sigma(W,\xi)]\cdot[\pi^{-1}(\overline{Z^-_{I\lambda}})]=[\Sigma(W,\xi)']\cdot [\pi^{-1}(\overline{Z^-_{I\lambda}})]
$$
is represented by a nonnegative cycle.
\qed

\section{A family of bases in which any Legendrian Thom polynomial
has positive expansion}\label{family}

We shall apply Theorem \ref{trans1} in the situation when all
$\alpha_i$ are equal to the same line bundle $\alpha$ (i.e.
$W=\alpha^{\oplus n}$) and $\alpha^{-m} \otimes \xi$ is globally
generated for $m\geq 3$.

\begin{example}\label{3cases} \rm
We shall consider the following three cases: the base is always
$X={\bf P}^n$ and

\begin{enumerate}
\item $$ \xi_1=\cO(-2)\,,\qquad \alpha_1=\cO(-1)\,,$$

\item $$ \xi_2=\cO(1)\,,\qquad \alpha_2={\bf 1}\,,$$

\item $$ \xi_3=\cO(-3)\,,\qquad \alpha_3=\cO(-1)\,,$$
\end{enumerate}
We obtain symplectic bundles $V_i=\alpha_i^{\oplus n}\oplus(\alpha_i^*\otimes\xi_i)^{\oplus n}$ with twisted symplectic forms $\omega_i$ for $i=1,2,3$.

These cases were for the authors the key examples supporting
their working conjecture that the assertion of Theorem \ref{main2}
holds true.

Case 1 was the subject of \cite[Remark 14]{MPW}, where the basis
related to the distinguished cell decomposition of
$LG(V_1,\omega_1)$ was investigated. In this case, for degrees $\le
n$, the cohomology $H^*(LG(V_1,\omega_1),\Z[\frac12])$ is isomorphic
to the ring of Legendrian characteristic classes tensored by
$\Z[\frac12]$.

In Case 2, the integral
cohomology $H^*(LG(V_2,\omega_2),\Z)$ is
isomorphic to the ring of Legendrian characteristic classes up to degree $n$.
The distinguished cell decomposition of $LG(V_2,\omega_2)$ gives us
another basis of cohomology.

In Case 3, the cohomology of $LG(V_3,\omega_3)$ is isomorphic, up to
degree $n$, to the ring of Legendrian characteristic classes,
provided we invert the number three this time.

The positivity property in Case 1 was known to us (see
\cite[Remark 14]{MPW}), whereas in Cases 2 and 3, it was Kazarian
who suggested the positivity. His conjecture was supported by
computation of all the Thom polynomials up to degree seven.

In general, $H^*(LG(V,\omega), {\Bbb Q})$ is isomorphic to the ring of Legendrian
characteristic classes up to degree $\dim W$ if $\xi$ is nontrivial. The case of $W=W^{(p,q)}$ and $\xi=\xi^{(p,q)}$, where
$$
\xi^{(p,q)}=\xi_2^{\otimes p}\otimes \xi_3^{\otimes q} \ \ \ \hbox{and} \ \ \ \alpha=\alpha^{(p,q)}=\alpha_2^{\otimes p}\otimes \alpha_3^{\otimes q}=\alpha_3^{\otimes q}
$$
($p,q$ are integers), will be used in Section \ref{parameter}.

\end{example}

\bigskip
To overlap all these three cases we consider the product
\begin{equation}\label{baza}
X:=\P^n\times \P^n
\end{equation}
and set
\begin{equation}\label{wiazka}
W:= p_1^* \cO(-1)^{\oplus n}\,,\qquad\xi:=p_1^*\cO(-3)\otimes p_2^*\cO(1)\,,
\end{equation}
where $p_i:X\to \P^n$, $i=1,2$, are the projections.  Restricting
the bundles $W$ and $\xi$ to the diagonal, or to the factors we
obtain the three cases considered above.
We should keep in mind that $X$ is an approximation of the classifying space $B(U(1)\times U(1))$. We   fix it because we apply  algebraic geometry methods.

\smallskip

The space $LG(V,\omega)$ has a distinguished cell decomposition
$Z^-_{I\lambda}$ where $I$ runs over strict partitions contained
in $\rho$, and $\lambda=(a,b)$ with $a$ and $b$ natural numbers
smaller than or equal to $n$. The classes of closures of the cells
of this decomposition give a basis of the homology of
$LG(V,\omega)$. The dual basis of cohomology (in the sense of linear algebra) is denoted by
\begin{equation}
e_{I,a,b}=[\overline{Z^-_{I,a,b}}]^*\,.
\end{equation}
By reasons of geometry, it is clear that the basis
$\{e_{I,a,b} \}$ consists of the classes represented by the cycles
$\overline{Z^+_{I,a,b}}$.


Let $v_1$ and $v_2$ be the first Chern classes of $p_1^*(\cO(1))$ and $p_2^*(\cO(1))$.
By the definition of $Z^+_{I,a,b}$, in $H^*(LG(V,\omega), \Bbb Z)$, we have
\begin{equation}\label{factor}
e_{I,a,b}=e_{I,0,0} \ v_1^a v_2^b\,.
\end{equation}
Moreover, we have
\begin{equation}\label{IOO}
e_{I,0,0}=[\Omega_I(F^+_\bullet,\xi)]\,.
\end{equation}

\medskip

With $X$, $W$ and $\xi$ as in (\ref{baza}) and (\ref{wiazka}), we have:

\begin{theorem}\label{main2} Let $\Sigma$ be a Legendre singularity class.
Then $[\Sigma(W,\xi)]$ has nonnegative coefficients in the basis $\{e_{I,a,b}\}$.
\end{theorem}
\proof
Let
$$
\iota: LG(V,\omega)\to {\gJ}^k(W,\xi)
$$
be the zero section, and
 $$
\iota^*: H^*( {\gJ}^k(W,\xi), {\Bbb Z}) \to H^*(LG(V,\omega),
{\Bbb Z})
$$
be the induced map on cohomology. We write
\begin{equation}
\iota^*[\Sigma(W,\xi)]=:\sum_{I,a,b}\gamma_{I,a,b}[\overline{Z^+_{I,a,b}}]\,,
\end{equation}
where $\gamma_{I,a,b}$ are integers.
We claim that  the coefficients $\gamma_{I,a,b}$ are nonnegative.
These coefficients are equal to
$$
\iota^*[\Sigma(W,\xi)]\cdot[\overline{Z^-_{I,a,b}}] \qquad {\rm (intersection\;in}\;LG(V,\omega)\,).$$
By the functoriality of the intersection product, the numbers $\gamma_{I,a,b}$ are equal to
 $$[\Sigma(W,\xi)]\cdot[\pi^{-1}(\overline{Z^-_{I,a,b}})] \qquad ({\rm intersection\;in}\;\gJ).$$
The vector bundle $\gJ$ on $LG(V,\omega)$ is equal to
$$
\tau^*\left(\bigoplus_{j=3}^{k+1}{\Sym}^j(W^*)\otimes
\xi\right)=\tau^*\left(\bigoplus_{j=3}^{k+1}{\Sym}^j({\bf 1}^n) \otimes
p_1^*\cO(j\moins 3)\otimes p_2^*\cO(1)\right).
$$
We see that $\gJ$ is globally generated. By Theorem \ref{trans1},
the desired intersections in $\gJ$ are nonnegative. \qed

\bigskip

The following computation will be needed later.

\begin{example}\label{exck}\rm  By \cite[Theorem 9.3]{PR} (see also \cite[Cor.~5]{KT})
if $\xi={\bf 1}$, $I=\{h\}$, then
$$[\Omega_h(F^+_\bullet,{\bf 1})]=c_h(R^*-F^+_{n+1-h})\,.$$
Hence
$$
[\Omega_h(F^+_\bullet,\xi)]=c_h(R^*\otimes\xi^{\frac12}-F^+_{
n+1-h}\otimes\xi^{-\frac 12})=c_h((R^*\otimes\xi-F^+_{
n+1-h})\otimes\xi^{-\frac 12})\,.
$$
Note that for any virtual bundle $E$ of dimension $h-1$ and for
any line bundle $\zeta$ we have $c_h(E\otimes\zeta)=c_h(E)$. Hence
$$[\Omega_h(F^+_\bullet,\xi)]=c_h(R^*\otimes\xi-F^+_{
n+1-h})\,.
$$
In our situation, ($W=\alpha^{\oplus n}$), the above formula can be written in the form
$$
[\Omega_h(F^+_\bullet,\xi)]=c_h(R^*\otimes\xi-W+\alpha^{\oplus
h-1})=c_h(W^*\otimes\xi-R+\alpha^{\oplus h-1})=c_h(A+\alpha^{\oplus h-1})\,.
$$
\end{example}

\section{The parameter $p/q$ and the basis for $p=q=1$}\label{parameter}

Fix a Legendre singularity class $\Sigma$.
By Theorem \ref{main2}, we know that the Thom polynomial of $\Sigma$,
evaluated at the Chern classes of
$$
A=W^*\otimes \xi-R
$$
and $c_1(\xi)=v_2-3v_1$, is a nonnegative $\Z$-linear combination of the following form:
$$
\cT^\Sigma =\sum_{I,a,b}\gamma_{I,a,b} \ e_{I,a,b}=
 \sum_{I,a,b}\gamma_{I,a,b} [\Omega_I(F^+_\bullet,\xi)]v_1^av_2^b\,.
$$

\smallskip

We want to find an additive basis of the ring of Legendrian characteristic classes
with the property that any Legendrian Thom polynomial is a nonnegative combination
of basis elements. To this end, we take a geometric model of the classifying space:
$LG(V^{(p,q)}, \omega^{(p,q)})$ (see the previous section), and the $Z^+$-basis which
is  dual to the $Z^-$-basis.
More precisely, dividing the cohomology ring $H^*(LG(V,\omega),\Qm)$ by the relation
\begin{equation}
q \cdot v_1=p \cdot v_2\,,
\end{equation}
that is specializing the parameters to $v_1=p \cdot t$,
$v_2=q \cdot t$, we obtain the ring $H^*(LG(V^{(p,q)}, \omega^{(p,q)}),{\Bbb Q})$,
isomorphic to the ring of Legendrian characteristic classes in degrees $\le n$ (provided
that $c_1(\xi)=v_2-3v_1$ is not specialized to 0.)

From Theorem \ref{main2}, we obtain:

\begin{theorem}\label{mainfamily}
If $p$ and $q$ are nonnegative, $q-3p\not=0$  then the Thom polynomial is a nonnegative combination of
the $[\Omega_I(F^+_\bullet,\xi)]\, t^i$'s.
\end{theorem}

The family $[\Omega_I(F^+_\bullet,\xi)]\, t^i$ is a one-parameter
family of bases depending on the parameter $p/q$.

\smallskip

Though the main theme of the present paper is the existence of
one-parameter family of bases in which every Legendrian Thom
polynomial has positive expansion, we shall give also some results
on algebraic form of $[\Omega_I(F^+_\bullet,\xi)]$.

First, we come back to Case 1 from Example \ref{3cases}.
This corresponds to fixing the parameter to be 1, i.e.
$p=1$ and $q=1$. This corresponds to setting $v_1=v_2=t$.
Geometrically, this means
that we study the restriction of the bundles $W$ and $\xi$ to the
diagonal of  $\P^n\times \P^n$, or, we study $W_1=\cO(-1)^{\oplus
n}$ and $\xi=\xi_1=\cO(-2)$. Set $\zeta:=\cO(-1)$. Then
$\xi^{\frac12}=\zeta$. From (\ref{rel}), we have
\begin{equation}\label{rel1}
A^*\otimes \xi^{\frac 1 2} + A\otimes \xi^{-\frac 1 2}=0\,.
\end{equation}

In general, when the bundle $\xi$ admits a square root $\zeta$ it
is convenient to give another description of the space
$LG(V,\omega)$. Let us define
$$
W'=W\otimes\zeta^{-1} \ \ \ \hbox{and} \ \ \
V'=V\otimes\zeta^{-1}\,.
$$
Then $V'$ is equipped with a symplectic form with constant
coefficients, and
$$V'=W'\oplus W'^*\,.$$
In our case,
$$
X=\P^n, \ \ \  W={\cal O}(-1)^{\oplus n}, \ \ \ \hbox{and} \ \ \
\xi={\cal O}(-2)\,.
$$
Then $W'$ and $V'$ become trivial bundles:
$$
W'={\bf 1}^n\,,\qquad V'={\bf 1}^{2n}\,.
$$
We have
\begin{equation}\label{ident}
LG(V,\omega)=LG(\C^{2n})\times \P^n\,,
\end{equation}
where $LG(\C^{2n})=LG(\C^{2n},\omega)$ and $\omega$ is the
standard nondegenerate symplectic form on $\C^{2n}$.

In algebraic expressions, we shall use $\Q$-functions of \cite{PR}
and their geometric interpretation from \cite{P1}. The reader can find in \cite[Sect.3]{MPW} a
summary of their properties in the notation which will be also used here.

Let $R'$ denote the tautological bundle on $LG(\C^{2n})$. Under the
identification (\ref{ident}), $R'$ pulled back from $LG(\C^{2n})$
to $LG(V,\omega)$ is equal to $R\otimes \zeta^{-1}$. We thus have
\begin{equation}\label{LR}
A\otimes\zeta^{-1}=W^*\otimes \zeta-R\otimes
\zeta^{-1}=W'^*-R'={\bf 1}^n-R'=R'^*-{\bf 1}^n\,.
\end{equation}
The distinguished cell decomposition of
$LG(V,\omega)=LG(\C^{2n})\times \P^n$ is of the product form.
By the cohomological properties of $\Q$-functions (\cite{P1}),
the basis of cohomology consists of the following
functions
$$\Q_I(R'^*)\cdot t^j=\Q_I(R^*\otimes\xi^{\frac12})\cdot
c_1(\xi^{-\frac12})^j\,,$$
where $I$ runs over strict partitions in $\rho$.

In this way, we obtain

\begin{theorem}\label{Qt}
The Thom polynomial for a Legendre singularity class $\Sigma$ is a combination:
\begin{equation}\label{alphaIj}
{\cal T}^\Sigma = \sum_{j\ge 0}\sum_I \alpha_{I,j} \ \Q_I(A\otimes
{\xi}^{-\frac 1 2})\cdot t^j \,,
\end{equation}
where
$$
t={\frac 1 2}c_1(\xi^*)\in H^2(X,\Z[{\frac 1 2}])\,, 
$$
$I$ runs over strict partitions in $\rho$,
and $\alpha_{I,j}$ are nonnegative integers.
\end{theorem}

\begin{remark} \rm We get the result announced in \cite[Remark 14]{MPW},
where it should read  ``$t=\frac 1 2 c_1(\xi^*)$'', and where we used
the notation $L^*-{\bf 1}^{\dim L}$ for the virtual bundle $A$.
\end{remark}

\section{Upper bounds of coefficients}\label{bounds}

In this section, we shall translate the positivity result in Theorem \ref{main2} into
restrictions for the coefficients of Thom polynomials.

Let us fix the value of the parameter. It turns out that the nonnegativity
of coefficients of the bases for some other values of the parameter
can imply upper bounds of the coefficients in the basis for the given parameter.
We plan to discuss this more systematically elsewhere. Here we consider
the 3 bases in degree two from Cases 1,2 and 3 in Example \ref{3cases}.
We shall call them the first, second and third basis, respectively.

Let $s=c_1(\xi)$ and let us list the classes of degree two. The first basis of
the Legendrian characteristic classes consists of:
$$
c_2(A+\alpha_1)=a_2+\frac12s\, a_1\,, \ \ \
 -\frac12 c_1(\xi_1)c_1(A)=-\frac12 s\,a_1 \,, \ \ \
 \left({-\frac12}c_1(\xi_1)\right)^2=\frac14s^2 \,,
$$
by Example \ref{exck}, since $c_1(\xi_1)=-2c_1(\cO(1))$ and
$c_1(\alpha_1)=-c_1(\cO(1))$.
The second basis is
$$
c_2(A)=a_2\,, \ \ \ s\, c_1(A)=s\, a_1\,, \ \ \ s^2\,
$$
since here $\alpha_2=1$.
The third basis is
$$
c_2(A+\alpha_3)=a_2+\frac13s\, a_1\,, \ \ \
 -\frac13 c_1(\xi_3)c_1(A)=-\frac13 s\,a_1 \,, \ \ \
 \left({-\frac13}c_1(\xi_3)\right)^2=\frac19s^2 \,,
$$
 since $c_1(\xi_3)=-3c_1(\cO(1))$ and $c_1(\alpha_3)=-c_1(\cO(1))$.

The Thom polynomial of the singularity  $A_3$ is of the form
$$
3(a_2+\frac12s\, a_1)-\kappa\,\frac12 s\,a_1
$$
with $\kappa\ge 0$, by the positivity in the first basis.
The positivity in the second basis gives the condition
$$
\frac32-\frac12\kappa\ge 0
$$
that is $\kappa\leq 3$.  When we write the Thom polynomial in the
third basis we see that the coefficient of $-\frac13 s\,a_1$ is
equal to $\kappa-1$. It follows that $\kappa\geq 1$.

Recall that the only Legendrian Thom polynomial of degree 2 is the
one of the singularity $A_3$, displayed in the first basis as:
$$
{\cal T}^{A_3}=3\Q_2+t \Q_1\,,
$$
i.e. with $\kappa=1$.

\medskip

Another upper bound for the coefficients of the expansions of Legendrian Thom
polynomials can be obtained by the method of Example \ref{hak}.

\section{Legendrian vs. classical Thom polynomials}\label{vs}

In this section, we shall use the basis from Case 1 in Example
\ref{3cases}, i.e., we put $t=v_1=v_2$.

\begin{proposition}\label{nonzero} For a nonempty stable Legendre singularity class $\Sigma$ the
Lagrangian Thom polynomial (i.e. ${\cal T}^\Sigma$ evaluated at $t=0$) is nonzero.
\end{proposition}
{\bf Proof.}
For a Legendre singularity class $\Sigma$ consider the associated singularity class
of maps $f:M \to C$ from $n$-dimensional manifolds to curves (see \cite[p.729]{Ka2} and \cite[p.123]{Kahab}).
We denote the related Thom polynomial by ${\it Tp}^\Sigma$.

According to \cite[pp.~708--709]{Ka2}, we have
\begin{equation}\label{fact}
{\it Tp}^\Sigma ={\cal T}^\Sigma \cdot c_n(T^*M\otimes f^*TC)\,.
\end{equation}
We know by \cite[Theorem 4]{PW} that the Thom polynomial ${\it Tp}^\Sigma$ is nonzero. Moreover,
it follows from the proof {\it (loc.cit.)} that ${\it Tp}^\Sigma$, specialized with
$f^*TC={\bf 1}$ i.e. $t=0$, is also nonzero. The assertion follows from the equation (\ref{fact}).
\qed

Consequently, we get an improvement of Theorem \ref{main2}.
\begin{corollary} For a nonempty stable Legendre singularity class $\Sigma$
the (Legendrian) Thom polynomial ${\cal T}^\Sigma$ is nonzero.
\end{corollary}

\begin{remark}\label{change} \rm Suppose that we are in the situation of (\ref{fact}).
We shall use the expansions of Legendrian Thom polynomials in the
basis for $v_1=v_2=t$, studied in the previous section.
Set $\xi:=f^*TC$. If $A=T^*M\otimes \xi-TM$, having the Thom
polynomial presented as in (\ref{alphaIj}), we want to compute the
$\Q$-functions of
$$
A\otimes \xi^{-\frac12} =(T^*M\otimes \xi-TM)^*\otimes
\xi^{-\frac12}=E^*-E\,,
$$
where $E=TM\otimes \xi^{-\frac 12}$.
Since for every strict partition $I$,
\begin{equation}
\Q_I(E^*-E)=Q_I(E^*)\,,
\end{equation}
where $Q_I$ denotes the classical {\it Schur $Q$-function}
\cite{Sch}, we get the desired expression by changing any
$\Q_I(A\otimes \xi^{-\frac12})$ to $Q_I(E^*)$.
\end{remark}

The following procedure mimics the passing from the LHS to RHS in
Eq.~(\ref{fact}), where ${\cal T}^\Sigma$ is given as a
${{\Z}[t]}$-combination of the $\Q_I(A\otimes \xi^{-\frac12})$'s,
and $Tp^\Sigma$ is to be written as a ${\Z}$-combination of the
Schur functions $S_J(T^*M-\xi^*)$ (for Schur functions in virtual
bundles we refer to \cite{L}, and for Schur $Q$-functions to
\cite{P1}).

\begin{procedure}\label{proc} \rm We start from a ${{\Z}[t]}$-combination of the polynomials
discussed in Remark \ref{change}: $\Q_I(A\otimes
\xi^{-\frac12})=Q_I(T^*M\otimes \xi^{\frac12})$.

\begin{itemize}

\item We write $Q_I(T^*M\otimes \xi^{\frac12})$ as a combination of the $S_J(T^*M\otimes \xi^{\frac12})$'s
(here we use a combinatorial rule from \cite{St} decomposing Schur $Q$-functions into $S$-functions);

\item we expand any $S_J(T^*M\otimes \xi^{\frac12})$ as a combination of the $S_K(T^*M)t^i$'s
(here we use a formula from \cite{L0} for decomposition of the Schur polynomials of twisted bundles);

\item we multiply the obtained combination by $c_n(T^*M\otimes \xi)$ (here we use the factorization
formula for Schur functions from \cite{BR}); we eventually get a combination of the $S_L(T^*M-\xi^*)$'s
with the coefficients being polynomials in $n$.

\end{itemize}

\end{procedure}

\begin{example} \rm \label{hak} We shall examine now how positivity of Schur function expansions
of Thom polynomials for mappings ${\C}^n \to {\C}$, proved in \cite{PW}, implies some upper bounds
on the coefficients of a Legendrian Thom polynomial in the expansion (\ref{alphaIj}).

Let us consider a degree 2 cohomology class already considered in Section \ref{bounds} of the form
\begin{equation}\label{QQ}
3\Q_2+\kappa t \Q_1\,,
\end{equation}
where $\kappa$ is an integer.

\smallskip

We fix $n\ge 2$. We apply Procedure \ref{proc} to (\ref{QQ}):
$$
\aligned &(3\Q_2+\kappa t \Q_1)(T^*M\otimes \xi^{\frac 1 2})\cdot
c_n(T^*M\otimes \xi)= \cr &(6(S_2+S_{1^2})+2\kappa t
S_1)(T^*M\otimes \xi^{\frac 1 2})\cdot c_n(T^*M\otimes \xi)= \cr
&\bigl(6(S_2+S_{1^2})+2t(\kappa-n)S_1+t^2n(n-2\kappa)
\bigr)(T^*M)\cdot c_n(T^*M\otimes \xi) \cr
\endaligned
$$
By the factorization formula, the last expression is equal to
\begin{equation}\label{schur}
6(S_{1^{n-1}3}+S_{1^{n-2}2^2})+(6n-\kappa)S_{1^n2}+(\frac32
n^2-\kappa \frac n2)S_{1^{n+2}}
\end{equation}
evaluated at $T^*M-\xi^*$.

Suppose that for each $n\ge 2$ we have in (\ref{schur}) a nonnegative combination of Schur functions,
i.e.
$$
6n-\kappa \ge 0 \ \ \ \hbox{and} \ \ \ \frac32 n^2-\kappa \frac
n2\ge 0 \,.
$$
This implies that $\kappa \le 6$.
\end{example}

\section{Examples of Legendrian Thom polynomials}\label{examples}

The Thom polynomials expanded in the basis $\{e_{I,a,b}\}$ \
(see Section \ref{family}) are (the summands in bold represent
the Lagrangian Thom polynomials):

\medskip

 \noindent $\bf A_2$: $\bf\Q_1$

 \noindent $\bf A_3$: ${\bf 3\Q_2}+v_2\Q_1$

\noindent $\bf A_4$: ${\bf
12\Q_3+3\Q_{21}}+(3v_1+7v_2)\Q_2+(v_1v_2+v_2^2)\Q_1$

\noindent $\bf D_4$: $\bf\Q_{21}$.
\medskip

\noindent The last equation means that the Thom polynomial of the
singularity $D_4$ written in all bases from the family is equal to
$\bf\Q_{21}$.

Similarly the Thom polynomial of the singularity $P_8$
in all bases from the family is equal to $\bf\Q_{321}$.
Next we have

\medskip

\noindent $\bf A_5$: ${\bf     60\Q_4+27\Q_{31}}+
      (6v_1+16v_2)\Q_{21}+
           (39v_1+47v_2)\Q_3+$

               $ (6v_1^2+22v_1v_2+12v_2^2)\Q_2+
(2v_1^2v_2+3v_1v_2^2+v_2^3)\Q_1$

\noindent $\bf D_5$: ${\bf 6\Q_{31}}+4v_2\Q_{21}$,

\medskip

\noindent and analogously to $D_5$,

\medskip

\noindent $\bf P_9$: ${\bf 12\Q_{421}}+12v_2\Q_{321}$.

\medskip

Let us specialize $v_1=v_2=t$.  The Thom polynomials for
singularities of codimensions lower or equal to six are listed in
\cite{MPW}. Here are the Legendrian Thom polynomials for the
consecutive singularities $A_8, D_8, E_8, X_9, P_9$ displayed in
the notation from (\ref{alphaIj}). Again, the summands in bold
represent the Lagrangian Thom polynomials. These formulas were
communicated to us by Kazarian.

$$
\aligned
&{\bf A_8:} \hskip0.95\hsize\cr
&{\bf 18840\Q_{61}+20160\Q_{7}+3123\Q_{421}+5556\Q_{43}+15564\Q_{52}}+ \cr
&t(71856\Q_{6}+3999\Q_{321}+55672\Q_{51}+34780\Q_{42})+ \cr
&t^2(64524\Q_{41}+24616\Q_{32}+105496\Q_{5})+t^3(36048\Q_{31}+81544\Q_{4})+ \cr
&t^4(8876\Q_{21}+34936\Q_{3})+t^5 7848\Q_{2}+t^6720\Q_1\,; \cr
\endaligned
$$
$$
\aligned
&{\bf D_8:} \hskip0.95\hsize\cr
&{\bf 1080\Q_{61}+315\Q_{421}+468\Q_{43}+1332\Q_{52}}+ \cr
&t(2754\Q_{42}+2952\Q_{51}+405\Q_{321})+t^2(1802\Q_{32}+3162\Q_{41})+ \cr
&t^31618\Q_{31}+t^4344\Q_{21}\,; \cr
\endaligned
$$
$$
\aligned
&{\bf E_8 :} \hskip0.95\hsize\cr
&{\bf 93\Q_{421}+108\Q_{43}+204\Q_{52}+72\Q_{61}}+t(99\Q_{321}+216\Q_{51}+414\Q_{42})+ \cr
&t^2(246\Q_{41}+246\Q_{32})+t^3126\Q_{31}+t^424\Q_{21}\,; \cr
\endaligned
$$
$$
\aligned
&{\bf X_9:} \hskip0.95\hsize\cr
&{\bf 18\Q_{52}+27\Q_{43}}+t(42\Q_{42}+6\Q_{51})+t^2(21\Q_{32}+11\Q_{41})+t^36\Q_{31}+t^4\Q_{21}\,; \cr
\endaligned
$$
\s
$$
\aligned
&{\bf P_9:} \hskip0.95\hsize\cr
&{\bf 12\Q_{421}}+t12\Q_{321}\,.\cr
\endaligned
$$

\s

The lower codimensional classes are displayed in \cite{MPW}, where the expression
for ${\bf A_7}$ should read:
$$
\aligned
&{\bf A_7:} \hskip0.95\hsize\cr
&{\bf 135 \Q_{321}+ 1275 \Q_{42}+ 2004\Q_{51}+  2520 \Q_{6}}+ \cr
&t(7092\Q_{5}+ 4439 \Q_{41}+ 1713 \Q_{32})+t^2(3545\Q_{31}+ 7868\Q_{4})+ \cr
&t^3(1106 \Q_{21}+4292\Q_{3})+t^4 1148\Q_{2}+t^5 120\Q_{1}\,. \cr
\endaligned
$$

\section{Appendix: A positivity result for Lagrange Grassmann bundle}

We shall give here a certain new positivity result which is a by-product
of our investigations in Section \ref{pos}. We follow the setup of
that section.

We assume here that
$X$ is homogeneous. For any automorphism of $X$ which is covered
by a map of $\xi$ and $\alpha_i$'s, we obtain an automorphism of
$LG(V,\omega)\to X$ transforming the fibers to fibers.

Inspired by the paper \cite{An} of Anderson, we consider some
``large group'' action. Assume that the line bundles:
$$
\alpha_i^*\otimes \alpha_j \ \ \hbox{for} \ i<j \ \ \ \hbox{and} \
\ \ \alpha_i^*\otimes\alpha_j^*\otimes\xi  \ \ \hbox{for all} \
i,j
$$
are globally generated. Consider the group $\Gamma B^-$ of global
sections of the bundle $B^- \to X$.

\begin{lemma} For each point $x\in X$, the restriction map from
$\Gamma B^-$ to the fiber $B^-_x$ is surjective.
\end{lemma}
\proof The group $B^-_x$ is generated by two subgroups:
\begin{itemize}
 \item $B_W^-$: the automorphisms of $W$ inducing automorphisms of
 $W^*\otimes\xi$ which preserve the flag $F^-_{\bullet}$,

\item $N^-$: the maps $W\to W^*\otimes \xi$ which belong to
$$
{\Sym}^2(W^*)\otimes\xi \subset W^*\otimes W^*\otimes
\xi=Hom(W,W^*\otimes \xi)\,.
$$
\end{itemize}
We identify the elements of $B_W ^-$ with matrix $\{b_{ij}\}$,
whose entries over $x\in X$ belong to the fiber
$$
Hom(\alpha_i,\alpha_j)_x=(\alpha_i^*\otimes\alpha_j)_x
$$
for $i\le j$, or are zero for $i>j$. The group bundle $B_W^-$ is
generated by the global sections. Similarly, the group bundle
$N^-$ is a quotient bundle of
$$
Hom(W,W^*\otimes\xi)\cong W^*\otimes W^*\otimes\xi
$$
isomorphic to ${\Sym}^2(W^*)\otimes\xi$; therefore it is globally
generated. This proves the lemma.
\qed

For a strict partition $J\subset \rho$, let us denote by $\Omega_J^-$
the space of the stratum $\Omega_J(F^-_\bullet,\xi) \to X$.

From the lemma, it follows

\begin{corollary} The group $\Gamma B^-$  acts on
$LG(V,\omega)$, preserving fibers, and in each fiber its orbits
coincide with the strata of the stratification $\{\Omega_J^- \}$.
\end{corollary}

\smallskip

Assume that $X$ is homogeneous with respect to a linear group $G$
and the transformation group acts on the line bundles $\xi$ and
$\alpha_i$. For instance, ~$X$ is a product of projective spaces,
and each line bundle involved is a tensor product of the
$\cO(j)$'s.

We define $H$ to be the subgroup of ${\rm Aut}(LG(V,\omega))$
generated by $\Gamma B^-$ and $G$ (it is the semidirect product
of these groups). The variety $H$ is irreducible.

From the above, it follows

\begin{lemma}\label{H}
The group $H$ acts transitively on each stratum
$\Omega^-_J$: $G$ transports any fiber to any other fiber,
and $\Gamma B^-$ acts transitively inside the fibers.
\end{lemma}

We now state the following positivity result.

\begin{theorem}\label{trans}
The intersection of any nonnegative cycle on $LG(V,\omega)$ with any
$\overline{Z^+_{I\lambda}}$ is represented by a nonnegative cycle.
\end{theorem}
\proof Let $Y$ be a nonegative cycle on $LG(V,\omega)$. We shall
find a translate $h\cdot Y$ by an element $h\in H$ which is
transverse to any $Z^+_{I\lambda}$.

By Lemma \ref{H}, we can use the Bertini-Kleiman transversality theorem \cite{Kl} for
$H$ acting on $\Omega^-_J$. By this theorem, there exists an open, dense subset
$U_{JI\lambda}\subset H$ with the following property:
if $h\in U_{JI\lambda}$ then $h\cdot (Y\cap \Omega^-_J)$ is transverse
to $Z^+_{I\lambda}\cap \Omega^-_J$. Set
\begin{equation}
U_J:=\bigcap_{I,\lambda} U_{JI\lambda}.
\end{equation}
We get an open, dense subset $U_J\subset H$ with the following property:
if $h\in U_J$, then $h\cdot (Y\cap \Omega^-_J)$ is transverse to any
$Z^+_{I\lambda}\cap \Omega^-_J$ (transversality in $\Omega^-_J$).
Since $\Omega^-_J$ is transverse to all strata $Z^+_{I\lambda}$ of
$LG(V,\omega)$, this transversality holds also in the whole ambient
space.
Set
\begin{equation}
U:=\bigcap_{{\rm strict} \ J \subset \rho} U_J.
\end{equation}
Pick $h\in U$. Then $Y'=h\cdot Y$ is transverse to all the $Z^+_{I\lambda}$'s.

The theorem now follows since by \cite[Sect. 8.2]{Fu} the intersection
$$
Y' \cdot [\overline{Z^+_{I\lambda}}]
$$
is represented by a nonnegative cycle.
\qed

\bigskip

\noindent
{\bf Acknowgledments.} We thank Maxim Kazarian for stimulating comments on the
first version of the present paper, and sending us his notes \cite{Ka3}.

The second author thanks the Tata Institute of Mathematical Sciences
in Mumbai and Institute of Mathematical Sciences in Chennai for hospitality.

The third author thanks the Institute of Mathematics of Polish Academy of Sciences
for hospitality.

\end{document}